\documentclass[a4paper,12pt]{article}
\usepackage{amsmath,amsthm,amssymb,latexsym,epic}
\usepackage{graphicx,enumerate}

\newtheorem{theorem}{Theorem}
\newtheorem{lemma}[theorem]{Lemma}

\newtheorem{remark}[theorem]{Remark}
\usepackage[all]{xy}

\newcommand{\IS}{{\mathcal {IS}}}
\newcommand{\dom}{{\mathrm {dom}}}

\renewcommand{\H}{{\mathcal {H}}}

\begin{document}
\title{Combinatorial Gelfand models for some semigroups 
and $q$-rook monoid algebras}
\author{Ganna Kudryavtseva and Volodymyr Mazorchuk}
\date{}

\maketitle

\begin{abstract}
Inspired by the results of \cite{APR}, we propose combinatorial 
Gelfand models for semigroup algebras of some finite semigroups, 
which include the symmetric inverse semigroup, the dual symmetric 
inverse semigroup, the maximal factorizable subsemigroup in the
dual symmetric inverse semigroup, and the factor power of the 
symmetric group. Furthermore we extend the Gelfand model for
the semigroup algebras of the symmetric inverse semigroup to a 
Gelfand model for the $q$-rook monoid algebra.
\end{abstract}

\section{Introduction}\label{s1}

Let $A$ be a finite-dimensional unital associative algebra over 
$\mathbb{C}$ and $M$ be an $A$-module. The module $M$ is said to be
a {\em Gelfand model} for $A$ if it is isomorphic to a multiplicity-free
direct sum of all simple $A$-modules. The present paper is inspired
by the results of \cite{APR}, where beautiful combinatorial Gelfand
models for the group algebra $\mathbb{C}S_n$ of the symmetric group, and 
for the  corresponding Hecke algebra $\mathbf{H}_n(q)$ are constructed.
We refer the reader to the list of references in \cite{APR} for the
history of the question and an account on known Gelfand models. 

The aim of the present paper is to extend the results of \cite{APR}
to some classes of finite semigroups, which include several inverse
generalizations of the symmetric group, in particular, the full
symmetric inverse semigroup $\mathcal{IS}_n$; and to the Hecke algebra
analogue for $\mathcal{IS}_n$, known as the $q$-rook monoid
algebra. The latter one has been recently defined by Solomon, \cite{So}, 
however, a special case already appeared in \cite{So1}. The $q$-rook 
monoid algebra has been studied by several authors, see  
\cite{DHP,HR,Ha,Pa,So} and references therein. Our motivation 
comes from an attempt to better understand the connection between 
the combinatorial and representation theoretical properties of these 
objects.

The paper is organized as follows: in Section~\ref{s2} we recall the
combinatorial Gelfand model for $\mathbb{C}S_n$, constructed in
\cite{APR}. In Section~\ref{s3} we show how the latter model can be used
to construct combinatorial Gelfand models for semigroup algebras
of those finite semigroups, for which each trace of a regular $\mathcal{D}$-class is an inverse semigroup in which maximal subgroups 
are direct sums of symmetric groups. Examples of such semigroups include
the symmetric inverse semigroup (see \cite[2.5]{GM2}), the dual 
symmetric inverse semigroup (see \cite{FL}), and the maximal factorizable
subsemigroup in the dual symmetric inverse  semigroup (see \cite{FL}).
Another, rather surprising, natural example is the factor power of the 
symmetric group (see \cite{GM0}), which, in particular, is not even regular.
In Section~\ref{s4} we recall (an appropriate modification of) 
the combinatorial Gelfand model for the Hecke algebra $\mathbf{H}_n(q)$,
constructed in \cite{APR}. Finally, in Section~\ref{s5} we extend the 
latter model to a combinatorial Gelfand model for the $q$-rook monoid 
algebra $\mathbf{I}_n(q)$ from \cite{So}. For $\mathbf{I}_n(q)$ we use 
the presentation from \cite{Ha}, which is different from the one used 
in \cite{So}.
\vspace{0.5cm}

\noindent
{\bf Acknowledgment.} This work was done during the visit of the
first author to Uppsala University, which was supported by the 
Royal Swedish Academy of Sciences  (KVA) and the Swedish Foundation for
International Cooperation in Research and Higher Education (STINT). For the
second author the research was partially supported by the Swedish Research
Council. The financial support of the Swedish Research Council, KVA and 
STINT and the hospitality of Uppsala University are gratefully acknowledged.
We thank Rowena Paget for her remarks on the original version of the
paper, especially on Section~\ref{s2}. We also thank the referee for very
useful remarks.

\section{Combinatorial Gelfand model for $\mathbb{C}S_n$}\label{s2}

Let $S_n$ be the symmetric group on $\{1,2,\dots,n\}$ and
$\mathcal{I}_n$ be the set of all involutions in $S_n$ (recall that
$\pi\in S_n$ is an involution provided that $\pi^2=\mathrm{id}$, in 
particular, the identity element $\mathrm{id}$ itself is an involution). 
For $\pi\in S_n$ we define the {\em inversion set} of $\pi$ as follows:
\begin{displaymath}
\mathrm{Inv}(\pi)=\{(i,j)\,:\, i<j \text{ and }\pi(i)>\pi(j)\}.
\end{displaymath}
For $w\in \mathcal{I}_n$ set
\begin{displaymath}
\mathrm{Pair}(w)=\{(i,j)\,:\, i<j \text{ and }w(i)=j\}.
\end{displaymath}
Set $\mathrm{Inv}_w(\pi)= \mathrm{Inv}(\pi)\cap \mathrm{Pair}(w)$
and $\mathrm{inv}_w(\pi)=|\mathrm{Inv}_w(\pi)|$. Finally, let $V_n$ 
be the vector space with the basis $\{\mathtt{I}_w:w\in \mathcal{I}_n\}$.

\begin{theorem}[\cite{APR}]\label{theorem1}
The assignment
\begin{displaymath}
\pi\cdot \mathtt{I}_w= (-1)^{\mathrm{inv}_w(\pi)}\mathtt{I}_{\pi w\pi^{-1}},
\quad\quad \pi\in S_n, w\in \mathcal{I}_n,
\end{displaymath}
defines on $V_n$ the structure of a $\mathbb{C}S_n$-module. Moreover,
this module is a Gelfand model for $\mathbb{C}S_n$.
\end{theorem}

\begin{remark}\label{remark1}
{\rm
Theorem~\ref{theorem1} extends to direct sums of symmetric groups in
a straightforward way.
}
\end{remark}

\begin{remark}\label{remark2}
{\rm
For $k=0,1,\dots,\lfloor\frac{n}{2}\rfloor$ let  $\mathcal{I}_n^k$ denote 
the subset of $\mathcal{I}_n$ consisting of all involutions, which can 
be written as a product of exactly $k$ pairwise different and commuting
transpositions. We obviously have that $\mathcal{I}_n$ is a disjoint union 
of the $\mathcal{I}_n^k$-s. Moreover, the linear span $V_n^k$ of 
$\{\mathtt{I}_w:w\in \mathcal{I}_n^k\}$ is invariant under the 
$\mathbb{C}S_n$-action for every $k$. 
The Robinson-Schensted correspondence, \cite[3.1]{Sa}, assigns to each
$\pi\in S_n$ a pair $(a(\pi),b(\pi))$ of standard Young tableaux of the
same shape. Moreover, $\pi\in S_n$ is an involution if and only if
$a(\pi)=b(\pi)$, \cite[Theorem~3.6.6]{Sa}. Using the properties of 
Viennot's shadow diagrams (see \cite[3.6]{Sa}) one shows that
two elements $w,w'\in \mathcal{I}_n$ belong to the same $\mathcal{I}_n^k$
provided that  $a(w)$ and $a(w')$ have the same shape. Using the main 
result of \cite{IRS} and tensoring with the sign representation one 
further shows that $V_n^k$ is isomorphic to the 
direct sum of Specht modules $S^{\lambda}$, where $\lambda$ runs 
through the set of all shapes of $a(w)$ for $w\in \mathcal{I}_n^k$.
}
\end{remark}

\section{Combinatorial Gelfand models for semigroup algebras
of some finite semigroups}\label{s3}

We use \cite{Hi} or \cite{GM2} as general reference for standard
notions from semigroup theory. Let $S$ be a finite semigroup and
$E(S)$ its set of idempotents. For $e\in E(S)$ consider the
$\mathcal{D}$-class $\mathrm{D}_e$ containing $e$. Then 
$\mathrm{D}_e\cup\{0\}$ with multiplication given by
\begin{displaymath}
a\star b=
\begin{cases} 
ab,& ab\in \mathrm{D}_e;\\
0, & \text{otherwise}, 
\end{cases}
\end{displaymath}
is called the {\em trace} of $\mathrm{D}_e$. From now on we 
assume that for every $e\in E(S)$
\begin{itemize}
\item the trace $\mathrm{D}_e\cup\{0\}$ is an  inverse semigroup,
\item the maximal subgroup $G_e$ of $S$, corresponding 
to  $e$, is a direct sum of symmetric groups.
\end{itemize}
Let $e_1,\dots,e_k\in E(S)$ be a fixed collection of idempotents, one 
for each $\mathcal{D}$-class. Let further $m_i$, $i=1,\dots,k$, denote 
the number of $\mathcal{L}$-classes inside $\mathrm{D}_{e_i}$. For 
each $i=1,\dots,k$ we fix an isomorphism of the group $G_{e_i}$ with  $S_{n^{(i)}_1}\oplus\dots\oplus S_{n^{(i)}_{l_i}}$ and 
an isomorphism $\varphi_i$ of $\mathrm{D}_{e_i}\cup\{0\}$ 
with the Brandt semigroup associated with the group 
$S_{n^{(i)}_1}\oplus\dots\oplus S_{n^{(i)}_{l_i}}$ and 
the cardinality $m_i$ (for details see \cite[\S 3.3]{CP}). 
This means that we have $\varphi_i(0)=0$ and
for any $x\in \mathrm{D}_{e_i}$ we have $\varphi_i(x)=(a,y,b)$, where $y\in 
S_{n^{(i)}_1}\oplus\dots\oplus S_{n^{(i)}_{l_i}}$, and
$a,b\in\{1,\dots,m_i\}$. The multiplication in the 
Brandt semigroup is given by
\begin{displaymath}
(a,y,b)\star (a',y',b')=
\begin{cases}
(a,yy',b'), & b=a';\\
0, & \text{otherwise}.
\end{cases}
\end{displaymath}
We set $\overline{\varphi}_i(x)=y$.

An element $w\in S$ will be called an {\em involution} provided that
$w\in G_e$ for some $e\in E(S)$ and $w^2=e$. Let $\mathcal{I}_S$
denote the set of all involutions in $S$, and $V_S$ the vector space 
with the basis $\{\mathtt{I}_w:w\in \mathcal{I}_S\}$. Our first result
is the following:

\begin{theorem}\label{theorem10}
Let $x\in S$ and $w\in \mathcal{I}_S$ be such that $w\in G_e$ for
some $e\in E(S)$ and $e\in\mathrm{D}_{e_i}$ for some
$i\in\{1,2,\dots,k\}$. Then  the assignment
\begin{equation}\label{eq:action}
x\cdot \mathtt{I}_w=
\begin{cases}
(-1)^{\mathrm{inv}_{\overline{\varphi}_i(w)}(\overline{\varphi}_i(xe))}
\mathtt{I}_{(xe) w(xe)^{-1}},&
xe\in \mathrm{D}_e;\\
0, & \text{ otherwise;}
\end{cases}
\end{equation}
defines on $V_S$ the structure of a $\mathbb{C}S$-module. Moreover,
this module is a Gelfand model for $\mathbb{C}S$.
\end{theorem}

\begin{proof}
Let $x,y\in S$ be such that $x\cdot (y\cdot \mathtt{I}_w) = 0$. We show 
that in this case we also have $xy\cdot \mathtt{I}_w = 0$. Suppose first 
that $y\cdot \mathtt{I}_w = 0$. Let $\leq_{\mathcal{J}}$ denote the 
natural partial order on $S$ associated with Green's $\mathcal{J}$-relation.
Then $ye<_{\mathcal{J}}e$ and thus $xye\leq_{\mathcal{J}}ye<_{\mathcal{J}}e$,
which means that $xy\cdot \mathtt{I}_w = 0$. Suppose now that 
$y\cdot \mathtt{I}_w \neq 0$, but 
$x\cdot (y\cdot \mathtt{I}_w) = 0$. Denote $v=(ye)w(ye)^{-1}$. Let $f$ be 
an idempotent such that $v\in G_f$. Then $x\cdot \mathtt{I}_v=0$, which 
implies that $xf<_{\mathcal{J}} f$. The definition of $f$ yields that 
$fye \mathcal{R} ye$, which, in turn, gives us that $fye=ye$. Thus we have
$$
xye=xfye\leq_{\mathcal{J}}xf<_{\mathcal{J}} f.
$$
It follows that $xy\cdot \mathtt{I}_w = 0$ as well.

Let now $x,y\in S$ be such that $x\cdot (y\cdot \mathtt{I}_w) \neq 0$. 
Denote $G=G_{e_i}$ and set $\overline{x}=\overline{\varphi}_i(x)$ 
for every $x\in \mathrm{D}_{e_i}$. Without loss of generality we will assume
that $\varphi_i(G_e)=(1,G,1)$ and that $\overline{g}=g$ whenever $g\in G_e$.
Since $y\cdot \mathtt{I}_w \neq 0$ it follows that $ye{\mathcal{J}}e$ and 
thus $ye{\mathcal{L}}e$. In particular, $\varphi_i(ye)=(l, \overline{ye}, 1)$
for some $l$. Then we have $\varphi_i((ye)^{-1})=(1,\overline{ye}^{-1},l)$, 
where $\overline{ye}^{-1}$ is the inverse of $\overline{ye}$ in $G$. 
Let $f=(ye)(ye)^{-1}$ and $v=(ye)w(ye)^{-1}$. Then 
$\varphi_i(f)=(l,e_i,l)$ and 
$\varphi_i(v)=(l,\overline{ye}w\overline{ye}^{-1},l)$.  Since $x\cdot 
\mathtt{I}_v \neq 0$ it follows that $xf{\mathcal{J}}f$ and thus 
$xf{\mathcal{L}}f$. Set $u=(xf)(ye)w(ye)^{-1}(xf)^{-1}$. Applying the 
equality $xfye=xye$ one shows that 
$\varphi_i(u)=(k,\overline{xye}w\overline{xye}^{-1},k)$ for some $k$.

Now the first part of the proof amounts to checking the equality 
\begin{equation}\label{eq1}
(-1)^{\mathrm{inv}_w(\overline{ye})}\cdot (-1)^{\mathrm{inv}_{\overline{ye}w\overline{ye}^{-1}}(\overline{xf})} 
=(-1)^{\mathrm{inv}_w(\overline{xye})}.
\end{equation}

Let $(i,j)\in\mathrm{Pair}(w)$. Suppose that
$(i,j)\in\mathrm{Inv}_w(\overline{ye})$, that is $i<j$ and 
$\overline{ye}(i)>\overline{ye}(j)$. If $\overline{xf}(\overline{ye}(i))<\overline{xf}(\overline{ye}(j))$ then we have that $(\overline{ye}(j),\overline{ye}(i))$ belongs to 
$\mathrm{Inv}_{\overline{ye}w\overline{ye}^{-1}}(\overline{xf})$, 
and at the same time $(i,j)\not\in \mathrm{Inv}_w(\overline{xye})$. 
If $\overline{xf}(\overline{ye}(i))>\overline{xf}(\overline{ye}(j))$ then
$(\overline{ye}(j),\overline{ye}(i))$ does not belong to 
$\mathrm{Inv}_{\overline{ye}w\overline{ye}^{-1}}(\overline{xf})$, 
and at the same time $(i,j)\in \mathrm{Inv}_w(\overline{xye})$.
Analogously one considers the case $(i,j)\not\in
\mathrm{Inv}_w(\overline{ye})$, and \eqref{eq1} follows.
Therefore $V_S$ is indeed a $\mathbb{C}S$-module.

We are left to show that $V_S$ is a Gelfand model for $S$. We will use 
the fact that simple modules over the complex semigroup algebra of a 
finite semigroup $S$ are in bijective correspondence with simple
modules of $G_{e_i}$, $1\leq i\leq k$ (see \cite[Chapter~5]{CP} or
\cite[Theorem~7]{GMS} for a more modern approach). In view 
of this, it is enough to show that for a maximal subgroup $G$ of $S$ and 
a simple ${\mathbb C}G$-module $V$ the corresponding ${\mathbb C}S$-module 
$V'$ is isomorphic to a submodule of $V_S$, and then to make sure that the 
sum of the dimensions of all $V'$-s equals the dimension of $V_S$.

Let $1\leq i\leq k$ and ${\mathcal I}_{n,i}$ be the set of involutions
 contained in maximal subgroups of the ${\mathcal D}$-class 
$\mathrm{D}_{e_i}$. Then the  linear span $V_S^i$ of all $\mathtt{I}_w$, 
$w\in {\mathcal I}_{n,i}$, is a direct summand of $V_S$. The dimension 
of this direct summand equals 
$m_i\cdot |\{\mathtt{I}_w: w\in\mathcal{I}_{G_{e_i}}\}|$. 
The action of $G_{e_i}$ on the linear span $V(G_{e_i})$ of 
$\{\mathtt{I}_w\,:\,w\in\mathcal{I}_{G_{e_i}}\}$, coinsides with 
the action from \cite{APR}, and thus from \cite[Theorem 1.1.2]{APR} and 
Remark~\ref{remark1} it follows that $V(G_{e_i})$ is a multiplicity-free
direct sum of all simple  $G_{e_i}$-modules. Let $V$ be a simple
direct summand of $V(G_{e_i})$ (as a $G_{e_i}$-module).
Suppose that the image of $G_{e_i}$ under $\varphi_i$ is $(1,G_{e_i},1)$ 
thus identifying $V$ with $(1,V,1)$ (the latter is  a subalgebra of the
semigroup algebra of the Brandt semigroup we work with). Then 
the vector space $\hat{V}=\oplus_{k=1}^{m_i} (k,V,1)$ is a simple $S$-module, 
corresponding to $V$, and by construction is a direct summand of $V_S^i$.
We have
\begin{displaymath}
\dim(\hat{V}) =m_i\cdot \dim (V)
\end{displaymath}
and hence $V_S^i$ is isomorphic to the multiplicity-free direct sum of 
all $\hat{V}$, where $V$ runs 
through the set of all simple $G_{e_i}$-modules. This completes the proof.
\end{proof}

Theorem~\ref{theorem10} applies to many semigroups.
Here are some examples:

\begin{itemize}
\item The symmetric inverse semigroup $\mathcal{IS}_n$ of all partial
injections on $\{1,2,\dots,n\}$ (also called the rook monoid), see 
\cite[2.5]{GM2}. The conditions are satisfied because of 
\cite[2.6 and 5.1]{GM2}.
\item The dual symmetric inverse semigroup $\mathcal{I}^*_n$ 
(or the monoid of block bijections) from \cite{FL}. The conditions are 
satisfied because of \cite[Theorem 2.2]{FL}.
\item The maximal factorizable submonoid of $\mathcal{I}^*_n$ 
(or the monoid of uniform block bijections) from
\cite{FL}. The conditions are satisfied because of \cite[Section 3]{FL}.
\item The factor power $\mathcal{FP}^+(S_n)$ from \cite{GM0,GM}.
Unlike the previous examples, this semigroup in {\em not} inverse, 
moreover, it is not even regular. However, all the required conditions are 
satisfied because of \cite[Theorem 1]{GM} and \cite{Ma}.
\end{itemize}

\section{Combinatorial Gelfand model for the Hecke algebra}\label{s4}

For a permutation $\pi\in S_n$ define the {\em support} of $\pi$ as follows:
\begin{displaymath}
{\mathrm{supp}}(\pi)=\{x\in\{1,2,\dots,n\}\,:\,\pi(x)\neq x\}.
\end{displaymath}
For $1\leq i< n$ let $s_i$ denote the transposition $(i,i+1)$.

For $q\in\mathbb{C}^*$ consider the {\em Hecke algebra} $\mathbf{H}_n(q)$,
which is a ${\mathbb C}$-algebra with generators $\{T_i: 1\leq i< n\}$ 
and defining relations
\begin{eqnarray}
(T_i-q)(T_i+1)&=&0, \quad\quad\quad\quad\quad 1\leq i< n; \label{h1}\\
T_iT_j&=&T_jT_i, \quad\quad\quad\,\,\,\,  1\leq i<j-1< n-1;\label{h2}\\
T_iT_{i+1}T_{i}&=&T_{i+1}T_iT_{i+1}, \quad  1\leq i< n-1.\label{h3}
\end{eqnarray}
We have $\mathbf{H}_n(1)\cong \mathbb{C}S_n$ canonically.
Let $V_{n,q}$ denote the formal linear span of 
$\{\mathtt{I}_w\,:\, w\in\mathcal{I}_n\}$.
The following theorem is a slightly modified version of 
\cite[Theorem 1.2.2]{APR}, which is more adjusted to our purposes. 

\begin{theorem}\label{thm201}
Let $1\leq i<n$ and $w\in S_n$ be an involution. The assignment
\begin{equation}\label{eq:action_H_n_q}
T_i\cdot \mathtt{I}_w=
\begin{cases}
q\mathtt{I}_w, &i,i+1\not\in {\mathrm{supp}}(w);\\
-\mathtt{I}_w, &i,i+1\in {\mathrm{supp}}(w);\\
\mathtt{I}_{s_iws_i}, &i\in {\mathrm{supp}}(w), 
i+1\not\in {\mathrm{supp}}(w);\\
q\mathtt{I}_{s_iws_i}+(q-1)\mathtt{I}_w, 
&i\not\in {\mathrm{supp}}(w), i+1\in {\mathrm{supp}}(w)
\end{cases}
\end{equation}
defines on $V_{n,q}$ the structure of an $\mathbf{H}_n(q)$-module.
If we additionally assume that $q$ is not a root of unity, then 
$V_{n,q}$ is a Gelfand model for $\mathbf{H}_n(q)$.
\end{theorem}

\section{Combinatorial Gelfand model for the $q$-rook monoid 
algebra}\label{s5}

For $q\in\mathbb{C}^*$ the {\em $q$-rook monoid algebra} $\mathbf{I}_n(q)$ 
is defined (see \cite{Ha}) as a ${\mathbb C}$-algebra with generators 
$\{T_i: 1\leq i< n\}\cup\{P_i: 1\leq i\leq n\}$  and 
defining relations \eqref{h1},  \eqref{h2}, 
\eqref{h3} and
\begin{eqnarray}
T_iP_j&=&P_jT_i=qP_j,\quad\quad\quad\quad\, 1\leq i<j\leq n;\label{r1}\\
T_iP_j&=&P_jT_i,\quad\quad\quad\quad\quad\quad\quad 
1\leq j<i\leq n-1;\label{r2}\\
P_i^{2}&=&P_i,\quad\quad\quad\quad\quad\quad\quad\quad 
1\leq i\leq n; \label{r3}\\
P_{i+1}&=&P_iT_iP_i-(q-1)P_i,\quad 1\leq i\leq n-1.\label{r4}
\end{eqnarray}

This is not a semigroup algebra of some ``$q$-rook monoid'' but 
rather  a one-parameter ($q$-) deformation of the semigroup algebra 
of the rook monoid,  in particular, $\mathbf{I}_n(1)\cong 
{\mathbb C}{\mathcal{IS}}_n$ canonically.
For generic $q$ there is a (non-canonical) isomorphism
$\mathbf{I}_n(q)\cong {\mathbb C}{\mathcal{IS}}_n$ (see \cite{So1,So,Pa}).

To proceed, we need to fix some notation. Let $X=\{1,2,\dots, n\}$.
${\mathcal{IS}}_n$ acts on $X$ in the standard way by partial permutations.
For a partial transformation $\alpha$ we denote by $\mathrm{dom}(\alpha)$
the {\em domain} of $\alpha$. 
For a subset $A$ of $X$ denote by $e_A$ the identity transformation of $A$, 
and by $G(A)$ the $\H$-class of $e_A$, which consists of all 
$\pi\in\IS_n$ whose domains and images are equal to $A$. 

Set $V=V_{q}^{\mathcal{IS}_n}$ to be the vector space with the basis
$\mathtt{I}_w$, where $w$ is an involution of $\mathcal{IS}_n$.

Let $A\subset X$. For every $\pi\in G(A)$ define $\psi_A(\pi)\in S_n$ 
as the element whose action on $A$ coincides with that of $\pi$, and 
which acts identically on the set $X\setminus A$. The map $\psi_A$ 
gives rise to a monomorphism $\overline{\psi}_A$ from the 
linear span of $\{\mathtt{I}_w\,:\,w\in\mathcal{I}_{G(A)}\}$ to 
$V_{n,q}$ defined on the basis via $\overline{\psi}_A(\mathtt{I}_w)=\mathtt{I}_{\psi(w)}$. 
If $i,i+1\in A$ and $w\in G(A)$ is  an involution, we define
$$
T_i\circ \mathtt{I}_w=\overline{\psi}_A^{-1}(T_i\cdot \overline{\psi}_A(\mathtt{I}_w)),
$$
where the action $T_i\cdot \mathtt{I}_{\psi_A(w)}$ is given by~\eqref{eq:action_H_n_q}.

Now, for every generator $T_i$, $1\leq i\leq n-1$, and $P_i$, 
$1\leq i\leq n$, of $\mathbf{I}_n(q)$ we define a linear transformation
of $V$ as follows:
\begin{equation}\label{eq:action_of_T}
T_i\cdot \mathtt{I}_w=
\begin{cases}
T_i\circ \mathtt{I}_w,& i,i+1\in \dom(w);\\
q\mathtt{I}_w, & i,i+1 \not\in \dom(w);\\
\mathtt{I}_{s_iws_i}, & i\in \dom(w), i+1\not\in \dom(w);\\
q\mathtt{I}_{s_iws_i}+(q-1)\mathtt{I}_w, & i\not\in \dom(w), i+1\in \dom(w);
\end{cases}
\end{equation}

\begin{equation}\label{eq:action_of_P}
P_i\cdot \mathtt{I}_w =
\begin{cases}
\mathtt{I}_w, & \dom(w)\subset \{i+1,\dots, n\};\\
0, & \dom(w)\not\subset \{i+1, \dots, n\}.
\end{cases}
\end{equation}

\begin{theorem}\label{th:gelfand_model_i_n_q}
The assignments \eqref{eq:action_of_T} and \eqref{eq:action_of_P} define on 
$V$ the structure of an $\mathbf{I}_n(q)$-module. If we additionally 
assume that $q$ is not a root of unity, then $V$ is 
a Gelfand model for $\mathbf{I}_n(q)$.
\end{theorem}

To prove the theorem we need some preparation. The group $S_n$ acts on 
$\mathcal{I}_{\IS_n}$ by conjugation. This action gives rise to an 
action of $S_n$ on $V$ defined as follows: $\pi\mathtt{I}_w\pi^{-1}=\mathtt{I}_{\pi w\pi^{-1}}$,
$w\in \mathcal{I}_{\IS_n}$. We will need the following technical lemma:

\begin{lemma}\label{lem:aux}
Let $w\in \mathcal{I}_{\IS_n}$.   
\begin{enumerate}[(a)]
\item \label{a1} If $i,i+1\in \dom(w)$, then $\pi (T_i\circ \mathtt{I}_w)\pi^{-1}=T_{i+1}\circ  \mathtt{I}_{\pi w\pi^{-1}}$ 
for any $\pi\in S_n$ such that $\pi(k)=k+1$ for $k=i,i+1$.
\item \label {a2} If $i,i+1\in\dom(w)$ and $|j-i|>1$, then 
$s_{j}(T_i\circ \mathtt{I}_w)s_{j}= T_{i}\circ \mathtt{I}_{s_{j}ws_{j}}$.
\end{enumerate}
\end{lemma}

\begin{proof}
Let $\psi=\psi_{\dom(w)}$ and $\tau=\psi_{\dom(\pi w\pi^{-1})}$.
Applying the definition of $\circ$ one reduces the  first equality to 
$$\pi(\overline{\psi}^{-1}(T_i\cdot\overline{\psi}(\mathtt{I}_{w})))
\pi^{-1}=\overline{\tau}^{-1}(T_{i+1}\cdot\overline{\tau}
(\mathtt{I}_{\pi w\pi^{-1}})).$$
Observe that for  $k=i,i+1$ we have $k\in \mathrm{supp}(\psi(w))$ if and only if $k+1\in \mathrm{supp}(\tau(\pi w\pi^{-1}))$. It follows that if $T_i\cdot(\overline{\psi}(\mathtt{I}_{w}))$ is a linear combination of some $\mathtt{I}_{\psi(u)}$'s then $T_{i+1}\cdot(\overline{\tau}(\mathtt{I}_{\pi w\pi^{-1}}))$ is the same linear combination of the corresponding $\mathtt{I}_{\tau(\pi u\pi^{-1})}$'s. Hence we are left to check the equality
$$\pi(\overline{\psi}^{-1}(\mathtt{I}_{\psi(u)}))\pi^{-1}=
\overline{\tau}^{-1}(\mathtt{I}_{\tau(\pi u\pi^{-1})})).$$
The latter equality reduces to $\pi\mathtt{I}_u\pi^{-1}=
\mathtt{I}_{\pi u\pi^{-1}}$,
which follows from the definitions. This proves \eqref{a1}.

To prove \eqref{a2} we set $\psi=\psi_{\dom(w)}$ and $\tau=\psi_{\dom(s_j ws_j)}$. The required equality reduces to
$$
s_j(\overline{\psi}^{-1}(T_i\cdot\overline{\psi}(\mathtt{I}_{w})))s_j=
\overline{\tau}^{-1}(T_{i}\cdot\overline{\tau}(\mathtt{I}_{s_jws_j})).
$$
Observe that for $k=i,i+1$ we have $k\in \mathrm{supp}(\psi(w))$ if and only if $k\in \mathrm{supp}(\tau(s_jws_j))$. It follows that if $T_i\cdot(\overline{\psi}(\mathtt{I}_{w}))$ is a linear combination of some $\mathtt{I}_{\psi(u)}$'s then $T_{i}\cdot(\overline{\tau}(\mathtt{I}_{s_jws_j}))$ is the same linear combination of the corresponding $\mathtt{I}_{\tau(s_jus_j)}$'s. Hence we are left to check the equality
$$s_j(\overline{\psi}^{-1}(\mathtt{I}_{\psi(u)}))s_j=
\overline{\tau}^{-1}(\mathtt{I}_{\tau(s_jus_j)})).$$
The latter equality reduces to $s_j\mathtt{I}_us_j=\mathtt{I}_{s_jus_j}$, 
which follows from the definitions. This completes the proof.
\end{proof}

\begin{proof}[Proof of Theorem \ref{th:gelfand_model_i_n_q}.]
First we show that $V$ is indeed an $\mathbf{I}_n(q)$-module. For this
we have to check the defining relations.

{\em Relation \eqref{h1}.} Observe that with respect to the fixed basis 
of $V$ the matrix corresponding to the action of $T_i$ is a direct sum 
of blocks of three possible types: $(q)$, $(-1)$, $\left(\begin{array}{ll}0&q\\ 1& q-1\end{array}\right)$. Each of these blocks satisfies \eqref{h1}.

{\em Relation \eqref{h2}.} There are sixteen possible cases depending on 
whether each of the elements $i,i+1,j,j+1$ belongs to $\dom(w)$ or not. 
The cases where $i,i+1\not\in\dom(w)$ or $j,j+1\not\in\dom(w)$ are trivial 
since then the action of $T_i$ or $T_{j}$ respectively just multiplies
the vectors $\mathtt{I}_w$, $\mathtt{I}_{s_iws_i}$, or $\mathtt{I}_{s_{j}ws_{j}}$
respectively by $q$. 

As $i$ and $j$ appear in \eqref{h2} symmetrically, we are left to consider 
six cases.

{\em Case 1.} $i,i+1,j,j+1\in \dom(w)$. We have
\begin{equation*}
\mathtt{I}_w \stackrel{T_i}{\mapsto} T_i\circ \mathtt{I}_w \stackrel{T_{j}}{\mapsto} T_j\circ (T_i\circ \mathtt{I}_w),\quad
\mathtt{I}_w \stackrel{T_j}{\mapsto} T_j\circ \mathtt{I}_w \stackrel{T_{i}}{\mapsto} T_i\circ (T_j\circ \mathtt{I}_w),
\end{equation*}
and the claim follows applying the relation
\eqref{h2} for $\mathbf{H}_n(q)$ (Theorem~\ref{thm201}).

{\em Case 2.} $i,i+1,j\in \dom(w)$, $j+1\not\in\dom(w)$.
We have
\begin{equation*}
\mathtt{I}_w \stackrel{T_i}{\mapsto}T_i\circ \mathtt{I}_w \stackrel{T_{j}}{\mapsto} s_j(T_i\circ \mathtt{I}_w)s_j,\quad
\mathtt{I}_w \stackrel{T_j}{\mapsto} \mathtt{I}_{s_jws_j} \stackrel{T_i}{\mapsto} T_i\circ\mathtt{I}_{s_jws_j},
\end{equation*}
and the claim follows from Lemma~\ref{lem:aux}\eqref{a2}.

{\em Case 3.} $i,i+1,j+1\in \dom(w)$, $j\not\in\dom(w)$.
We have
\begin{equation*}
\mathtt{I}_w \stackrel{T_i}{\mapsto}T_i\circ \mathtt{I}_w \stackrel{T_{j}}{\mapsto} qs_j(T_i\circ \mathtt{I}_w)s_j+(q-1)T_i\circ \mathtt{I}_w,
\end{equation*}
\begin{equation*}
\mathtt{I}_w \stackrel{T_j}{\mapsto} q\mathtt{I}_{s_jws_j}+(q-1)\mathtt{I}_w \stackrel{T_i}{\mapsto}
qT_i\circ(\mathtt{I}_{s_jws_j})+(q-1)T_i\circ \mathtt{I}_w,
\end{equation*}
and the claim follows from Lemma \ref{lem:aux}\eqref{a2}.

{\em Case 4.} $i,j\in \dom(w)$, $i+1,j+1\not\in\dom(w)$. We have
\begin{equation*}
\mathtt{I}_w \stackrel{T_i}{\mapsto} \mathtt{I}_{s_iws_i}\stackrel{T_j}{\mapsto} \mathtt{I}_{s_js_iws_is_j},\quad
\mathtt{I}_w \stackrel{T_j}{\mapsto} \mathtt{I}_{s_jws_j}\stackrel{T_i}{\mapsto} \mathtt{I}_{s_is_jws_js_i},
\end{equation*}
and the claim follows as $s_js_i=s_is_j$.

{\em Case 5.} $i,j+1\in \dom(w)$, $i+1,j\not\in\dom(w)$. We have
\begin{equation*}
\mathtt{I}_w \stackrel{T_i}{\mapsto} \mathtt{I}_{s_iws_i}\stackrel{T_j}{\mapsto} q\mathtt{I}_{s_js_iws_is_j}+(q-1)\mathtt{I}_{s_iws_i},
\end{equation*}
\begin{equation*}
\mathtt{I}_w \stackrel{T_j}{\mapsto} q\mathtt{I}_{s_jws_j}+(q-1)\mathtt{I}_w\stackrel{T_i}{\mapsto}
q\mathtt{I}_{s_is_jws_js_i}+(q-1)\mathtt{I}_{s_iws_i},
\end{equation*}
and the claim follows as $s_js_i=s_is_j$.

{\em Case 6.} $i+1,j+1\in \dom(w)$, $i,j\not\in\dom(w)$. We have
\begin{equation*}
\mathtt{I}_w \stackrel{T_i}{\mapsto} q\mathtt{I}_{s_iws_i}+(q-1)\mathtt{I}_w\stackrel{T_j}{\mapsto}
q^2\mathtt{I}_{s_js_iws_is_j}+q(q-1)\mathtt{I}_{s_iws_i}+q(q-1)
\mathtt{I}_{s_jws_j}+(q-1)^2\mathtt{I}_w,
\end{equation*}
\begin{equation*}
\mathtt{I}_w \stackrel{T_j}{\mapsto} q\mathtt{I}_{s_jws_j}+(q-1)\mathtt{I}_w\stackrel{T_i}{\mapsto}
q^2\mathtt{I}_{s_is_jws_js_i}+q(q-1)\mathtt{I}_{s_jws_j}+q(q-1)
\mathtt{I}_{s_iws_i}+(q-1)^2\mathtt{I}_w,
\end{equation*}
and the claim follows as $s_js_i=s_is_j$.

{\em Relation \eqref{h3}.}
We consider eight possible cases depending on whether the elements 
$i,i+1,i+2$ belong to $\dom(w)$ or not.

{\em Case 1.} $i,i+1,i+2\in \dom(w)$. This follows immedeately
from  Theorem~\ref{thm201}.

{\em Case 2.} $i,i+1\in \dom(w)$, $i+2\not\in \dom(w)$. We have
\begin{equation*}\mathtt{I}_w \stackrel{T_i}{\mapsto} T_i\circ \mathtt{I}_w \stackrel{T_{i+1}}{\mapsto}  s_{i+1}(T_i\circ \mathtt{I}_w)s_{i+1}\stackrel{T_i}{\mapsto} \\s_is_{i+1}(T_i\circ \mathtt{I}_w)s_{i+1}s_i,
\end{equation*}
\begin{equation*}
\mathtt{I}_w \stackrel{T_{i+1}}{\mapsto} \mathtt{I}_{s_{i+1}ws_{i+1}} \stackrel{T_i}{\mapsto}  \mathtt{I}_{s_is_{i+1}ws_{i+1}s_i} \stackrel{T_{i+1}}{\mapsto}\\
T_{i+1}\circ \mathtt{I}_{s_is_{i+1}ws_{i+1}s_i},
\end{equation*} 
and the claim follows applying Lemma \ref{lem:aux}\eqref{a1} 
for $\pi=s_is_{i+1}$.

{\em Case 3.} $i,i+2\in \dom(w)$, $i+1\not\in \dom(w)$. We have 
\begin{equation}\label{a3}
\mathtt{I}_w \stackrel{T_i}{\mapsto} \mathtt{I}_{s_iws_i} \stackrel{T_{i+1}}{\mapsto} T_{i+1}\circ \mathtt{I}_{s_iws_i} \stackrel{T_i}{\mapsto} 
qs_i(T_{i+1}\circ \mathtt{I}_{s_iws_i})s_i+ 
(q-1)T_{i+1}\circ \mathtt{I}_{s_iws_i},
\end{equation}
\begin{multline}\label{b3}
\mathtt{I}_w \stackrel{T_{i+1}}{\mapsto} q\mathtt{I}_{s_{i+1}ws_{i+1}}+(q-1)\mathtt{I}_w \stackrel{T_i}{\mapsto}
qT_i\circ \mathtt{I}_{s_{i+1}ws_{i+1}}+(q-1)\mathtt{I}_{s_iws_i} \stackrel{T_{i+1}}{\mapsto}\\
qs_{i+1}(T_i\circ \mathtt{I}_{s_{i+1}ws_{i+1}})s_{i+1}+
(q-1)T_{i+1}\circ \mathtt{I}_{s_iws_i}.
\end{multline}
The claim now follows from the equality
\begin{displaymath}
s_is_{i+1}(T_i\circ \mathtt{I}_{s_{i+1}ws_{i+1}})s_{i+1}s_i= 
T_{i+1}\circ \mathtt{I}_{s_iws_i},
\end{displaymath}
which, in turn, follows from Lemma~\ref{lem:aux}\eqref{a1}.

{\em Case 4.} $i+1,i+2\in \dom(w)$, $i\not\in \dom(w)$. We have 
\begin{multline}\label{a4}
\mathtt{I}_w \stackrel{T_i}{\mapsto} q\mathtt{I}_{s_iws_i}+(q-1)\mathtt{I}_w \stackrel{T_{i+1}}{\mapsto}
q^2\mathtt{I}_{s_{i+1}s_iws_is_{i+1}}+q(q-1)\mathtt{I}_{s_iws_i}+(q-1)T_{i+1}\circ \mathtt{I}_w\\ \stackrel{T_i}{\mapsto}
q^2T_i\circ \mathtt{I}_{s_{i+1}s_iws_is_{i+1}}+q(q-1)\mathtt{I}_w+(q-1)qs_i(T_{i+1}\circ \mathtt{I}_w)s_i+
(q-1)^2 T_{i+1}\circ \mathtt{I}_w \\ =: A+B+C+D,
\end{multline}
\begin{multline}\label{b4}
\mathtt{I}_w \stackrel{T_{i+1}}{\mapsto} T_{i+1}\circ \mathtt{I}_w \stackrel{T_i}{\mapsto} qs_i(T_{i+1}\circ \mathtt{I}_w)s_i+(q-1)T_{i+1}\circ \mathtt{I}_w \stackrel{T_{i+1}}{\mapsto}\\
q^2s_{i+1}s_i(T_{i+1}\circ  \mathtt{I}_w)s_is_{i+1}+q(q-1)s_i(T_{i+1}\circ \mathtt{I}_w)s_i+\\ +(q-1)T_{i+1}\circ(T_{i+1}\circ \mathtt{I}_w)=:E+F+G.
\end{multline}
We have $C=F$. Further, $B+D=G$ follows from \eqref{h1}. Finally,
$A=E$ follows from Lemma \ref{lem:aux}\eqref{a1} applying arguments
similar to those used in the previous case.

{\em Case 5.} $i+2\in \dom(w)$, $i,i+1\not\in \dom(w)$. We have
\begin{multline}\label{a5}
\mathtt{I}_w \stackrel{T_i}{\mapsto} q\mathtt{I}_w 
\stackrel{T_{i+1}}{\mapsto} q^2\mathtt{I}_{s_{i+1}ws_{i+1}}+
q(q-1)\mathtt{I}_w \stackrel{T_i}{\mapsto}\\
q^3\mathtt{I}_{s_is_{i+1}ws_{i+1}s_i}+q^2(q-1)
\mathtt{I}_{s_{i+1}ws_{i+1}}+q^2(q-1)\mathtt{I}_w,
\end{multline}
\begin{multline}\label{b5}
\mathtt{I}_w \stackrel{T_{i+1}}{\mapsto} q\mathtt{I}_{s_{i+1}ws_{i+1}}
+(q-1)\mathtt{I}_w \stackrel{T_i}{\mapsto}\\ 
q^2\mathtt{I}_{s_is_{i+1}ws_{i+1}s_i} +q(q-1)\mathtt{I}_{s_{i+1}ws_{i+1}} 
+q(q-1)\mathtt{I}_w \stackrel{T_{i+1}}{\mapsto}\\
q^3\mathtt{I}_{s_is_{i+1}ws_{i+1}s_i}+q(q-1)\mathtt{I}_w+
(q-1)q^2\mathtt{I}_{s_{i+1}ws_{i+1}}+(q-1)^2q\mathtt{I}_w.
\end{multline}
The right hand sides of \eqref{a5} and \eqref{b5} are equal.

{\em Case 6.} $i+1\in \dom(w)$, $i,i+2\not\in \dom(w)$. We have
\begin{multline*}\label{a6}
\mathtt{I}_w \stackrel{T_i}{\mapsto} q\mathtt{I}_{s_iws_i}+(q-1)\mathtt{I}_w \stackrel{T_{i+1}}{\mapsto} 
q^2\mathtt{I}_{s_iws_i}+(q-1)\mathtt{I}_{s_{i+1}ws_{i+1}} \stackrel{T_i}{\mapsto}\\
q^2\mathtt{I}_w + q(q-1)\mathtt{I}_{s_{i+1}ws_{i+1}};
\end{multline*}
\begin{equation*}\label{b6}
\mathtt{I}_w \stackrel{T_{i+1}}{\mapsto} \mathtt{I}_{s_{i+1}ws_{i+1}} \stackrel{T_i}{\mapsto} q\mathtt{I}_{s_{i+1}ws_{i+1}} \stackrel{T_{i+1}}{\mapsto} q^2\mathtt{I}_w + q(q-1)\mathtt{I}_{s_{i+1}ws_{i+1}}.
\end{equation*}

{\em Case 7.} $i\in \dom(w)$, $i+1,i+2\not\in \dom(w)$. In this case we have
\begin{equation*}\label{a7}
\mathtt{I}_w \stackrel{T_i}{\mapsto} \mathtt{I}_{s_{i}ws_{i}}\stackrel{T_{i+1}}{\mapsto} \mathtt{I}_{s_{i+1}s_iws_is_{i+1}} \stackrel{T_i}{\mapsto} q\mathtt{I}_{s_{i+1}s_iws_is_{i+1}};
\end{equation*}
\begin{equation*}\label{b7}
\mathtt{I}_w \stackrel{T_{i+1}}{\mapsto} q\mathtt{I}_w\stackrel{T_i}{\mapsto} q\mathtt{I}_{s_{i}ws_{i}}\stackrel{T_{i+1}}{\mapsto} q\mathtt{I}_{s_{i+1}s_iws_is_{i+1}}.
\end{equation*}

{\em Case 8.} $i,i+1,i+2\not\in \dom(w)$. 
In this case the actions of both $T_iT_{i+1}T_i$ and $T_{i+1}T_iT_{i+1}$ map $\mathtt{I}_w$ to $q^3\mathtt{I}_w$.

{\em Relation \eqref{r1}.} We consider two possible cases.

{\em Case 1.} $\dom(w)\subset\{j+1,\dots, n\}$. Then we also have $i,i+1\not\in\dom(w)$. Therefore, $T_i\cdot \mathtt{I}_w=q\mathtt{I}_w$ 
and thus  $T_iP_j\cdot\mathtt{I}_w=P_jT_i\cdot\mathtt{I}_w=
qP_j\cdot\mathtt{I}_w=q\mathtt{I}_w.$

{\em Case 2.} $\dom(w)\not\subset\{j+1,\dots, n\}$. In this case 
it follows from
the definitions that $T_iP_j\cdot\mathtt{I}_w=P_jT_i\cdot\mathtt{I}_w=
qP_j\cdot\mathtt{I}_w= 0$.

{\em Relation \eqref{r2}.}  We consider two possible cases.

{\em Case 1.} $\dom(w)\subset\{j+1,\dots, n\}$. Then
$T_iP_j\cdot\mathtt{I}_w=P_jT_i\cdot\mathtt{I}_w=T_i\cdot\mathtt{I}_w$
as $T_i\cdot \mathtt{I}_w$ is a linear combination of 
$\mathtt{I}_w$ and $\mathtt{I}_{s_iws_i}$ and both $\dom(w)$ and
$\dom(s_iws_i)$ are contained in $\{j+1,\dots, n\}$.

{\em Case 2.} $\dom(w)\not\subset\{j+1,\dots, n\}$. In this case
$T_iP_j\cdot\mathtt{I}_w=P_jT_i\cdot\mathtt{I}_w=0$ follows from the
definitions.

{\em Relation \eqref{r3}} follows immediately from the definitions.

{\em Relation \eqref{r4}.} Consider three possible cases.

{\em Case 1.} $\dom(w)\subset \{i+2,\dots, n\}$. In this case we also have $i,i+1\not\in\dom(w)$. Therefore, $T_i\cdot \mathtt{I}_w=q\mathtt{I}_w$ 
and hence 
$$P_{i+1}\cdot \mathtt{I}_w=(P_iT_iP_i-(q-1)P_i)\cdot \mathtt{I}_w=
\mathtt{I}_w.$$

{\em Case 2.} $\dom(w)\subset \{i+1,\dots, n\}$ and $i+1\in \dom(w)$. In this case $P_{i+1}\cdot\mathtt{I}_w=0$ and $P_{i}\cdot\mathtt{I}_w=\mathtt{I}_w$. We are left to check that $P_iT_i\cdot \mathtt{I}_w=(q-1)I_w$. Observe that
$P_i\cdot\mathtt{I}_{s_iws_i}=0$ as $\dom(s_iws_i)\not\subset
\{i+1,\dots, n\}$. Taking this into account we  obtain
$$
P_iT_i\cdot \mathtt{I}_w=P_i\cdot (q\mathtt{I}_{s_iws_i}+(q-1)\mathtt{I}_w)=(q-1)\mathtt{I}_w.
$$

{\em Case 3.} $\dom(w)\not\subset \{i+1,\dots, n\}$. This case
follows directly from the definitions. This completes the proof of 
the fact that $V$ is an ${\mathbf{I}}_n(q)$-module.

To prove that $V$ is a Gelfand model for ${\mathbf{I}}_n(q)$ we use
the results of \cite{Pa} and arguments analogous to those used in
the second part of the proof of Theorem~\ref{theorem10}. 
Set $J_0={\mathbf{I}}_n(q)$, $J_k= {\mathbf{I}}_n(q)P_k{\mathbf{I}}_n(q)$,
$k=1,\dots,n$, and  $J_{n+1}=0$. 
In \cite[Section~4]{Pa} it is shown that there is an algebra
isomorphism  
\begin{displaymath}
{\mathbf{I}}_n(q)\cong \bigoplus_{k=0}^n J_k/J_{k+1}
\end{displaymath}
and that $J_k/J_{k+1}$ is isomorphic to 
$\mathbf{M}_{l_k}(\mathbb{C})\otimes \mathbf{H}_{n-k}(q)$ (where
$\mathbf{M}_{l_k}(\mathbb{C})$ is the matrix algebra and $l_k=\binom{n}{k}$),
with a multiplication, which is appropiately twisted to take into account
powers of $q$, which may appear when multiplying in $J_k/J_{k+1}$. This
twist can be forgotten using \cite[Corollary~15]{Pa}.

For $k=0,\dots,n$ let $V^{(k)}$ denote the subspace of $V$, spanned
by all $\mathtt{I}_w$, where $|\dom(w)|=n-k$. Denote also by 
$\tilde{V}^{(k)}$ the subspace of $V^{(k)}$, spanned
by all $\mathtt{I}_w$, where $\dom(w)=\{1,2,\dots,n-k\}$.
Finally, we identify $\mathbf{H}_{n-k}(q)$ with the subalgebra of 
${\mathbf{I}}_n(q)$, generated by $T_i$, $i<n-k$. From the definitions 
we immediately obtain:
\begin{itemize}
\item $V=\oplus_{k=0}^n V^{(k)}$ as ${\mathbf{I}}_n(q)$-modules;
\item $J_{k+1}\cdot V^{(k)}=0$.
\end{itemize}
Furthermore, from the definitions we
have that $\tilde{V}^{(k)}$ is an $\mathbf{H}_{n-k}(q)$-module
and, moreover, from Theorem~\ref{thm201} we have that 
$\tilde{V}^{(k)}$ is a Gelfand model of $\mathbf{H}_{n-k}(q)$.

As $\dim V^{(k)}= l_k\cdot \dim \tilde{V}^{(k)}$, we obtain that 
$V^{(k)}$ is a Gelfand model 
for $\mathrm{M}_{l_k}(\mathbb{C})\otimes \mathbf{H}_{n-k}(q)$.
Adding everything up we thus get that $V$ is a Gelfand model for 
${\mathbf{I}}_n(q)$.  This completes the proof.
\end{proof}

\vspace{1cm}

\noindent
G.K.: Department of Mechanics and Mathematics, Kyiv Taras Shev\-chen\-ko
University, 64, Volodymyrska st., 01033, Kyiv, UKRAINE,
e-mail: {\em akudr\symbol{64}univ.kiev.ua}
\vspace{0.5cm} 

\noindent
V.M.: Department of Mathematics, Uppsala University, SE 471 06,
Uppsala, SWEDEN, e-mail: {\em mazor\symbol{64}math.uu.se}
\end{document}